\documentclass[12pt]{amsart}
\usepackage{amssymb}
\usepackage{dsfont}
\usepackage{amscd}
\usepackage{mathabx}
\usepackage[mathscr]{euscript} 
\usepackage[all]{xy}
\usepackage{url}
\usepackage{stmaryrd}
\usepackage{comment}
\usepackage[retainorgcmds]{IEEEtrantools}
\usepackage{hyperref}
\title{Positive characteristic Ax-Schanuel}
\author[P. KOWALSKI]{Piotr Kowalski$^{\spadesuit}$}
\thanks{$^{\spadesuit}$ Supported by the Narodowe Centrum Nauki grant no.  2021/43/B/ST1/00405 and by the T\"{u}bitak 1001 grant no. 119F397.}
\address{$^{\spadesuit}$Instytut Matematyczny\\
Uniwersytet Wroc{\l}awski\\
Wroc{\l}aw\\
Poland}
\email{pkowa@math.uni.wroc.pl} \urladdr{http://www.math.uni.wroc.pl/\textasciitilde pkowa/ }

 \DeclareMathOperator{\fr}{Fr}

\DeclareMathOperator{\td}{trdeg} 
 \DeclareMathOperator{\alg}{alg}

\DeclareMathOperator{\trd}{trdeg}
\DeclareMathOperator{\inv}{inv}

\DeclareMathOperator{\rat}{rat}

\DeclareMathOperator{\dcf}{DCF}

\newtheorem{theorem}{Theorem}[section]

\theoremstyle{definition}
\newtheorem{definition}[theorem]{Definition}
\newtheorem{example}[theorem]{Example}
\newtheorem{remark}[theorem]{Remark}

\begin{document} 
\newcommand{\lili}{\underleftarrow{\lim }}
\newcommand{\coco}{\underrightarrow{\lim }}
\newcommand{\twoc}[3]{ {#1} \choose {{#2}|{#3}}}
\newcommand{\thrc}[4]{ {#1} \choose {{#2}|{#3}|{#4}}}
\newcommand{\Zz}{{\mathds{Z}}}
\newcommand{\Ff}{{\mathds{F}}}
\newcommand{\Cc}{{\mathds{C}}}
\newcommand{\Rr}{{\mathds{R}}}
\newcommand{\Nn}{{\mathds{N}}}
\newcommand{\Qq}{{\mathds{Q}}}
\newcommand{\Kk}{{\mathds{K}}}
\newcommand{\Pp}{{\mathds{P}}}
\newcommand{\ddd}{\mathrm{d}}
\newcommand{\Aa}{\mathds{A}}
\newcommand{\dlog}{\mathrm{ld}}
\newcommand{\ga}{\mathbb{G}_{\rm{a}}}
\newcommand{\gm}{\mathbb{G}_{\rm{m}}}
\newcommand{\gaf}{\widehat{\mathbb{G}}_{\rm{a}}}
\newcommand{\gmf}{\widehat{\mathbb{G}}_{\rm{m}}}
\newcommand{\ka}{{\bf k}}
\newcommand{\ot}{\otimes}
\newcommand{\si}{\mbox{$\sigma$}}
\newcommand{\ks}{\mbox{$({\bf k},\sigma)$}}
\newcommand{\kg}{\mbox{${\bf k}[G]$}}
\newcommand{\ksg}{\mbox{$({\bf k}[G],\sigma)$}}
\newcommand{\ksgs}{\mbox{${\bf k}[G,\sigma_G]$}}
\newcommand{\cks}{\mbox{$\mathrm{Mod}_{({A},\sigma_A)}$}}
\newcommand{\ckg}{\mbox{$\mathrm{Mod}_{{\bf k}[G]}$}}
\newcommand{\cksg}{\mbox{$\mathrm{Mod}_{({A}[G],\sigma_A)}$}}
\newcommand{\cksgs}{\mbox{$\mathrm{Mod}_{({A}[G],\sigma_G)}$}}
\newcommand{\crats}{\mbox{$\mathrm{Mod}^{\rat}_{(\mathbf{G},\sigma_{\mathbf{G}})}$}}
\newcommand{\crat}{\mbox{$\mathrm{Mod}^{\rat}_{\mathbf{G}}$}}
\newcommand{\cratinv}{\mbox{$\mathrm{Mod}^{\rat}_{\mathbb{G}}$}}
\newcommand{\ra}{\longrightarrow}
\newcommand{\bdcf}{B-\dcf}
\makeatletter
\providecommand*{\cupdot}{%
  \mathbin{%
    \mathpalette\@cupdot{}%
  }%
}
\newcommand*{\@cupdot}[2]{%
  \ooalign{%
    $\m@th#1\cup$\cr
    \sbox0{$#1\cup$}%
    \dimen@=\ht0 %
    \sbox0{$\m@th#1\cdot$}%
    \advance\dimen@ by -\ht0 %
    \dimen@=.5\dimen@
    \hidewidth\raise\dimen@\box0\hidewidth
  }%
}

\providecommand*{\bigcupdot}{%
  \mathop{%
    \vphantom{\bigcup}%
    \mathpalette\@bigcupdot{}%
  }%
}
\newcommand*{\@bigcupdot}[2]{%
  \ooalign{%
    $\m@th#1\bigcup$\cr
    \sbox0{$#1\bigcup$}%
    \dimen@=\ht0 %
    \advance\dimen@ by -\dp0 %
    \sbox0{\scalebox{2}{$\m@th#1\cdot$}}%
    \advance\dimen@ by -\ht0 %
    \dimen@=.5\dimen@
    \hidewidth\raise\dimen@\box0\hidewidth
  }%
}
\makeatother

\def\Ind#1#2{#1\setbox0=\hbox{$#1x$}\kern\wd0\hbox to 0pt{\hss$#1\mid$\hss}
\lower.9\ht0\hbox to 0pt{\hss$#1\smile$\hss}\kern\wd0}

\def\ind{\mathop{\mathpalette\Ind{}}}

\def\notind#1#2{#1\setbox0=\hbox{$#1x$}\kern\wd0
\hbox to 0pt{\mathchardef\nn=12854\hss$#1\nn$\kern1.4\wd0\hss}
\hbox to 0pt{\hss$#1\mid$\hss}\lower.9\ht0 \hbox to 0pt{\hss$#1\smile$\hss}\kern\wd0}

\def\nind{\mathop{\mathpalette\notind{}}}


\maketitle

\begin{abstract}
This expository paper is written in celebration of Boris Zilber's 75th birthday. We discuss Ax-Schanuel type statements focusing on the case of positive characteristic.
\end{abstract}

\section{Introduction}
During the Spring 2005 Isaac Newton Institute Program ``Model Theory and Applications to Algebra and Analysis'' in Cambridge,
I learnt that I will be a MODNET postdoc with Boris Zilber in Oxford for the academic year 2005/2006. Still in Cambridge, Boris suggested to me to start thinking on ``positive characteristic versions of Ax's theorem''. In this expository paper, I will describe what has happened next.

It may be a good moment for a general disclaimer. This is an expository paper representing my experience with respect to the Boris' suggestion above and I do not claim that this paper describes adequately the state of the art in the vast area of Ax-Schanuel type problems. In particular, comparatively very little will be said about the amazing developments of Jonathan Pila (and many others) regarding the modular version of Ax-Schanuel and its applications to diophantine problems, most notable the Andr\'{e}-Oort conjecture. I will write more about it in Section \ref{secchar0}.

This paper is organized as follows. In Section \ref{secchar0}, we describe the history of this circle of problems in the case of characteristic $0$. In Section \ref{seccharp}, we focus on the positive characteristic case and present some of the results I obtained regarding to the Boris' suggestion above. In Section \ref{secspecul}, we speculate on some recent ideas regarding general forms of Ax-Schanuel and its Hasse-Schmidt differential versions.

I would like to thank Jakub Gogolok for his comments and to the referee for the very careful and useful report.

\section{Characteristic zero}\label{secchar0}

In this section, we summarize the characteristic $0$ situation regarding the Ax-Schanuel problems. The disclaimer from the Introduction above applies mostly here.
\subsection{Results}

In 1960's Schanuel formulated two conjectures \cite[pages 30–31]{lang1966}: one about transcendence of complex numbers \cite[(S)]{ax71} and one about transcendence of power series \cite[(SP)]{ax71}. We state them below.

\smallskip

{\bf Schanuel's Conjecture (Complex Numbers)}
\\
\textit{Let  $x_1,\ldots,x_n\in \Cc$ be linearly independent over $\Qq$. Then
 $$\trd_{\Qq}(x_1,\ldots,x_n,e^{x_1},\ldots,e^{x_n})\geqslant n.$$}

{\bf Schanuel's Conjecture (Power Series)}
\\
\textit{Let $x_1,\dots ,x_n\in X\Cc\llbracket X\rrbracket$ be linearly independent over $\Qq$. Then
$$\trd_{\Cc(X)}(x_1,\dots ,x_n,e^{x_1},\dots ,e^{x_n})\geqslant n.$$}

Complex Numbers Schanuel's Conjecture is open even for $n=2$, since (using Euler's identity $e^{i\pi}+1=0$) it covers the open problem of algebraic independence of $\pi$ and $e$ and it is even still unknown whether $\pi+e$ is irrational (it is phrased as a ``Candidate for the Most Embarrassing Transcendence Question in Characteristic Zero'' in \cite{Brow})! Power Series Schanuel's Conjecture was proved by Ax (\cite[(SP)]{ax71}).
\\
\\
Ax also showed in \cite{ax71} the following differential version of Power Series Schanuel's Conjecture, which was actually used by Ax to show the other statements from \cite{ax71}.

\vspace{3mm}

{\bf Differential Ax-Schanuel Theorem \cite[(SD)]{ax71}}
\\
\textit{Let $(K,\partial)$ be a differential field of characteristic 0 and $C$ be its field of constants.
For $x_1,\ldots,x_n\in K$ and $y_1,\ldots,y_n\in K^*$, if
$$\partial x_1=\frac{\partial y_1}{y_1},\ldots,\partial x_n=\frac{\partial y_n}{y_n}$$ and
$\partial x_1,\ldots,\partial x_n$ are $\Qq$-linearly independent, then
$$\trd_{C}(x_1,\ldots,x_n,y_1,\ldots,y_n)\geqslant n+1.$$}

\begin{remark}\label{passege}
There are the following passages between the power series and the differential version of Ax's theorem above.
\begin{enumerate}
  \item Since the ring of power series has a natural differential structure, the differential version implies the power series version.

  \item Going the other way is more subtle. Seidenberg's embedding theorem \cite{Seid} says that any finitely generated differential field of characteristic $0$ differentially embeds into the differential field of meromorphic functions on an open subset of $\Cc$. Using this theorem, one can reduce the differential version of Ax's theorem to the power series one (this is explained in detail around Theorem 4.1 in \cite{FrSc} and in Section 2.5 in \cite{PT}).
\end{enumerate}
These passages apply to the more complicated cases of analytic (or formal) Ax-Schanuel statements versus the differential ones as well. Such more complicated cases will be described below.
\end{remark}
In a subsequent paper written one year later \cite{ax72}, Ax proved the following general geometric result about the dimension of intersections of algebraic subvarieties of complex algebraic groups with analytic subgroups.

\smallskip

{\bf Ax's theorem on the dimension of intersections \cite[Theorem 1]{ax72}}
\\
\textit{Let $G$ be an algebraic group over the field of complex numbers $\Cc$. Let $\mathcal{A}$ be a complex analytic subgroup of $G(\Cc)$ and $V$ be an irreducible algebraic subvariety of $G$ over $\Cc$. We assume that $\mathcal{K}:=\mathcal{A}\cap V(\Cc)$ is Zariski dense in $V(\Cc)$. Then there is an analytic subgroup $\mathcal{B}\subseteq G(\Cc)$ containing $V(\Cc)$ and $\mathcal{A}$ such that
$$\dim(\mathcal{B})\leqslant \dim(\mathcal{A})+\dim(V)-\dim(\mathcal{K}).$$}
This theorem implies Power Series Schanuel's Conjecture by taking:
\begin{itemize}
  \item $G$ as the product of the vector group $\ga^n$ and the torus $\gm^n$,
  \item $\mathcal{A}$ as the $n$-th Cartesian power of the graph of the exponential map,
  \item $V$ as the algebraic locus of the tuple $(x_1,\dots ,x_n,e^{x_1},\dots ,e^{x_n})$.
\end{itemize}
Ax's theorem on the dimension of intersections applies also (more generally) to the case of the exponential map on a semi-abelian variety \cite[Theorem 3]{ax72}. The consequences of Ax's theorem on the dimension of intersections go beyond the case of the exponential map, for example this theorem applies to the case of analytic maps between the multiplicative group and an elliptic curve. We state it precisely below, since this statement is amenable for a possible transfer to the positive characteristic case (see Remark \ref{lookpos}).
\begin{theorem}\label{risingax}
Let
$$\gamma:\gm(\Cc)\to E(\Cc)$$
be an analytic epimorphism where $E$ is an elliptic curve. Let $x_1,\dots ,x_n\in 1+X\Cc\llbracket X\rrbracket$ be multiplicatively independent. Then
$$\trd_{\Cc(X)}\left(x_1,\dots ,x_n,\gamma(x_1),\dots ,\gamma(x_n)\right)\geqslant n.$$
\end{theorem}
After Ax's work in 1970's, Brownawell and Kubota \cite{BrKu} proved a version of the differential Ax's theorem in the
case of elliptic curves, and then Kirby \cite{ki} generalized it to arbitrary semi-abelian varieties. These
results were not included in \cite{ax72}, however they are closely related using the ``passages'' from Remark \ref{passege}. Bertrand  \cite{ber} extended \cite[Theorem 3]{ax72} to commutative algebraic groups not having vector quotients (e.g. a maximal non-split vectorial extensions of a semi-abelian variety).

Differential Ax's theorem \cite[(SD)]{ax71} is generalized further to “very
non-algebraic formal maps” in \cite[Theorem 5.5]{{K5}}. This generalization includes a differential version of Bertrand’s result and a differential Ax-Schanuel type result about raising to non-algebraic powers on an algebraic
torus \cite[Theorem 6.12]{K5}. We state it below in the power series case (see Remark \ref{passege}), since this statement has a positive characteristic interpretation (see Remark \ref{lookpos}). Before the statement, we note that for $x\in 1+X\Cc\llbracket X\rrbracket$ and $\alpha\in \Cc$, we define
$$x^{\alpha}:=\exp(\alpha \log(x)),$$
where $\exp,\log\in \Qq\llbracket X\rrbracket$ are the standard formal power series corresponding to the exponential and the logarithmic maps.
\begin{theorem}\label{risingax}
Suppose that $\alpha\in \Cc$ is such that $[\Qq(\alpha):\Qq]>n$, Let $x_1,\dots ,x_n\in 1+X\Cc\llbracket X\rrbracket$ be multiplicatively independent. Then
$$\trd_{\Cc(X)}\left(x_1,\dots ,x_n,x_1^{\alpha},\dots ,x_n^{\alpha}\right)\geqslant n.$$
\end{theorem}

We will describe now briefly modular analogues of Ax's theorem. Our disclaimer from the Introduction applies very much here. Ax-Schanuel statements may go beyond the context of group homomorphisms, the first example here is the $j$-function map:
$$j:\mathbb{H}\to \Cc,$$
where $\mathbb{H}$ is the upper half plane. The linear independence assumption from Ax's theorem is replaced with \emph{modular independence}. Pila's notes \cite{Pilanotes} contain an excellent comprehensive survey of the state of the art in this field up to Year 2013. Such results have very important diophantine applications such as:
\begin{itemize}
  \item another proof of the Manin-Mumford conjecture (\cite{pzmm});

  \item the first unconditional proof of the Andr\'{e}-Oort conjecture for $\Cc^n$ (\cite{Pila11});

  \item a recent proof of the full Andr\'{e}-Oort conjecture for Shimura varieties (\cite{PSTEG});
  \end{itemize}
Following a suggestion by the referee, we would like to point out that only the Ax-Lindemann-Weierstrass type of results are needed in Manin-Mumford and Andr\'{e}-Oort, while Ax-Schanuel (in fact, a weak form of it) is used in Zilber-Pink type problems.

In \cite{cfn1}, the Ax-Lindemann-Weierstrass theorem with derivatives for the uniformizing functions of genus zero Fuchsian groups of the first kind is shown.  This result is used in \cite{cfn1} to answer a question of Painlev\'{e} from 1895.
\begin{remark}\label{lookpos}
We analyze now which statements of the Ax-Schanuel results discussed above are transferrable to the positive characteristic case. We would like to mention that all the analytic Ax-Schanuel type results over $\Cc$ may be replaced with their formal counterparts over an arbitrary field $C$, which was already done by Ax: the reader is advised to compare Ax's theorem on the dimension of intersections with its formal counterpart \cite[Theorem 3]{ax72}, which will be stated in a more general form in Section \ref{seccharp}. Let us recall the set-up first.
\begin{definition}[Bochner \cite{Boc}]\label{deffor}
An $n$-dimensional \emph{formal group} (law) over $C$ is a tuple of power series $F\in
C\llbracket X,Y\rrbracket^{\times n}$ ($|X|=|Y|=|Z|=n$) satisfying:
\begin{itemize}
\item $F(0,X)=F(X,0)=X$,

\item $F(X,F(Y,Z))=F(F(X,Y),Z).$
\end{itemize}
A \emph{morphism} from an $n$-dimensional formal group $G$ into an $m$-dimensional formal group
$F$ is a tuple of power series $f\in C\llbracket X\rrbracket^{\times m}$ such that:
\begin{itemize}
\item $F(f(X),f(Y))=f(G(X,Y))$.
\end{itemize}
\end{definition}
\noindent
 There is a well-known formalization functor $G\mapsto \widehat{G}$ (see pages 5 and 13 in \cite{manin})
from the category of algebraic groups to the category of formal groups.

Such characteristic $0$ formal statements seem to be transferrable to the positive characteristic context in the cases when the corresponding formal maps exist.
\begin{enumerate}
  \item The very original version of Ax-Schanuel does not look transferable, since there are no reasonable exponential maps in positive characteristic (we will briefly touch on the Drinfeld context at the end of Section \ref{seccharp}).

  \item Therefore, other analytic maps need to be considered. ``Analytic'' may be replaced with ``formal'' (as mentioned above) and then the closest one to the exponential map which survives in the case of positive characteristic seems to be the formal isomorphism between the multiplicative group and an ordinary elliptic curve.

  \item The other types of such maps come from raising to powers in the multiplicative group.

\end{enumerate}
Items $(2)$ and $(3)$ above will be discussed in the positive characteristic context in Section \ref{seccharp}.
\end{remark}

\subsection{Motivations and applications}\label{22}
In \cite{Zil}, Zilber  used Differential Ax's Theorem
 to prove \emph{Weak CIT}, which is a weak version of \emph{Conjecture on Intersection with
Tori} (CIT), which was also stated in \cite{Zil}. CIT is a finiteness statement about intersections of subtori of
a given torus with certain subvarieties of this torus. Weak CIT was used in \cite{BaHiMPWa}
to produce a characteristic 0 \emph{bad field}. The existence of such a field was an open problem in model theory for almost 20 years.

Regarding the positive characteristic case, Weak CIT does not hold and Zilber formulated a conjectural statement in \cite{zilsur} (the very last statement of \cite{zilsur}).
It is still open whether a bad field in the positive characteristic case exists, however, Wagner showed \cite{Wa1} that its existence in the case of characteristic $p>0$ implies the existence of infinitely many \emph{$p$-Mersenne primes}, which is an open problem in number theory, but it is widely believed that there are finitely many of them (for each individual prime $p$). Therefore, the existence of bad fields in positive characteristic looks very unlikely. However, pursuing the following path of research still looks interesting:
\begin{enumerate}
  \item prove positive characteristic versions of Ax-Schanuel;

  \item show a version of Weak CIT in positive characteristic using Item $(1)$;

  \item construct a version of a bad field in positive characteristic using Item $(2)$;

  \item check the possible number-theoretic consequences of results obtained in Item $(3)$.
\end{enumerate}
As was mentioned in the previous subsection, Jonathan Pila and others used Ax-Schanuel type results to show different versions of the Andr\'{e}-Oort conjecture, see e.g. \cite{Pila11}, \cite{T}, \cite{cfn1}, and \cite{PSTEG}.

There are also model-theoretic consequences of results of Ax-Schanuel type and we would like to point out some of them.
\begin{itemize}
  \item In \cite{ki}, Kirby used his version of an Ax-Schanuel statement to obtain the complete first order theories of the exponential differential equations of semiabelian varieties which arise from an amalgamation construction in the style of Hrushovski.

  \item In \cite{aslan2}, Aslanyan did a version of the above for the $j$-function in place of the exponential function on semiabelian varieties.

  \item Freitag and Scanlon used Ax-Lindemann-Weierstrass to establish strong minimality and triviality of the differential equation of the $j$-function in \cite{FrSc}. This was generalized by Aslanyan in \cite{aslan1} to a more general and formal setting.

  \item In \cite{cfn1} and \cite{cfn3}, the authors go in a quite an opposite way: they first establish strong minimality using differential Galois theory, then use Zilber's trichotomy to prove triviality, then use that to establish Ax-Lindemann-Weierstrass and later Ax-Schanuel. That is, they give a new proof to Ax-Schanuel for the $j$-function and in fact for all Fuchsian automorphic functions.
\end{itemize}

\section{Positive characteristic}\label{seccharp}

The first (to my knowledge) positive characteristic Ax-Schanuel result concerns additive power series. Interestingly, it is not included in the cases considered in Remark \ref{lookpos}, because such formal maps have no counterpart in the characteristic 0 case, since any additive formal power series in characteristic zero is linear, so it is ``not interesting''. This positive characteristic additive Ax-Schanuel result is explained in detail below.

For any commutative algebraic group $G$, we have the following two (usually non-commutative) rings:
\begin{enumerate}
  \item  the ring of \emph{algebraic endomorphisms} (that is: endomorphisms of $G$ in the original category of algebraic groups), denoted $\mathrm{End}_{\mathrm{algebraic}}(G)$;
  \item  the ring of \emph{formal endomorphisms} (that is: endomorphisms of the formalization of $G$, as below Definition \ref{deffor}, in the category of formal groups), denoted $\mathrm{End}_{\mathrm{formal}}(G)$.
\end{enumerate}
Let $C$ be a field of characteristic $p>0$ and $\ga$ denote the additive group scheme over $C$. We consider the following two rings.
\begin{itemize}
  \item  The ring of additive polynomials over $C$ (with composition), which we denote by $C[\fr]$. This is also the  skew polynomial ring over $C$ and we have the following ring isomorphism:
      $$\mathrm{End}_{\mathrm{algebraic}}(\ga)\cong C[\fr].$$

  \item The ring of additive  power series over $C$ (with composition), which we denote by $C\llbracket \fr \rrbracket$. We have the following ring isomorphism:
      $$\mathrm{End}_{\mathrm{formal}}(\ga)\cong C\llbracket \fr \rrbracket.$$
\end{itemize}
These rings are commutative if and only if $C=\Ff_p$ and then they are also domains (isomorphic to the rings of polynomials or the ring of power series). We denote the fraction field of $\Ff_p[\fr]$ by $\Ff_p(\fr)$. We state below the main theorem of \cite{K6}.
\\
\\
{\bf Ax-Schanuel for additive power series \cite[Theorem 1.1]{K6}}
\\
\textit{Let $F$ be an additive power series over $\Ff_p$ and assume that
$$[\Ff_p(\fr)(F):\Ff_p(\fr)]>n.$$
 Let $x_1,\ldots,x_n\in t\Ff_p\llbracket t\rrbracket$ be linearly independent over $\Ff_p[\fr]$. Then we have:
$$\td_{\Ff_p}(x_1,\ldots,x_n,F(x_1),\ldots,F(x_n))\geqslant n+1.$$}
We will describe a general Ax-Schanuel result from \cite{K8}, which is valid in all characteristics. We need two technical assumptions. Before stating them, we will try to motivate them. One of the crucial properties (used in the proofs in \cite{ax72}) of analytic homomorphisms between algebraic groups is that they take invariant algebraic differential forms into invariant algebraic differential forms. The first technical assumption below, which is absolutely necessary, is both formalizing and generalizing this crucial property. Regarding the second assumption, the exponential map gives a formal isomorphism between any commutative algebraic (and even formal) group in the case of characteristic zero and a Cartesian power of the additive group. This is false in the positive characteristic case, for example there is no formal isomorphism between the additive and the multiplicative group (no exponential map in positive characteristic!). To make the proofs work, we still need to impose an additional assumption in the positive characteristic case, to mimic the above characteristic 0 situation. The 1-dimensional group $H$ in this assumption plays the role of $\ga$ and we need to put some extra conditions on $H$ which are true for $\ga$. We would prefer to avoid this second assumption, however, we were unable to do so in \cite{K8}.
\begin{enumerate}
  \item We define a \emph{special} formal map as one which ``resembles a homomorphism'' in the sense that it takes invariant differential forms into the ``usual'' differential forms (before taking the completion, see \cite[Def. 3.10]{K8}). In the positive characteristic case, the notion of differential forms has to be replaced by the Vojta's notion (\cite{voj}) of \emph{higher differential forms}, see \cite[Remark 5.18(3)]{K8}.

  \item We say that a commutative algebraic group $A$ defined over the field $C$ of characteristic $p$ is ``good''  (see \cite[Def. 3.4]{K8}), if there is a one-dimensional algebraic group $H$ over $C$ such that we have the following (in the case of $p=0$, we drop Item (c) below).
\begin{enumerate}
\item $\widehat{A}\cong \widehat{H^n}$.

\item The map $\mathrm{End}(\widehat{H})\to \mathrm{End}_C(\Omega^{\inv}_H)$($=C$) is onto.

\item $H$ is \emph{$\Ff_p$-isotrivial} i.e. $H\cong H^{\fr}$.
\end{enumerate}
\end{enumerate}
To motivate the next result and give a general feeling regarding ``what is it about'', we quote from \cite{K8} the following, where ``The main theorem of this paper'' refers to Theorem \ref{short}.
\\
\\
``A continuous map between Hausdorff spaces which is constant on a dense set is constant everywhere. The same principle applies to an algebraic map between algebraic varieties and to the Zariski topology (which is not Hausdorff). However, if we mix categories there is no reason for this principle to hold, e.g. there are non-constant \emph{analytic} maps between \emph{algebraic} varieties which are constant on a Zariski dense subset. The main theorem of this paper roughly says that the principle above can be saved for certain formal maps (resembling homomorphisms) between an algebraic variety and an algebraic group at the cost of replacing the range of the map with its quotient by a formal subgroup of the controlled dimension.''
\begin{theorem}\label{short}
Let $V$ be an algebraic variety, $\mathcal{K}$ a Zariski dense formal subvariety of $V$, $A$ a ``good'' commutative algebraic group and $\mathcal{F}:\widehat{V}\to \widehat{A}$ a special formal map. Assume $\mathcal{F}$ vanishes on $\mathcal{K}$. Then there is a formal subgroup $\mathcal{C}\leqslant \widehat{A}$ such that $\mathcal{F}(\widehat{V})\subseteq \mathcal{C}$ and
$$\dim(\mathcal{C})\leqslant \dim(V)-\dim(\mathcal{K}).$$
\end{theorem}
As a consequence of Theorem \ref{short}, we obtained in \cite{K8} a result which is parallel to Ax-Schanuel for additive power series, where an additive power series (that is: a formal endomorphism of the additive group) is replaced with a ``multiplicative'' power series (that is: a formal endomorphism of the multiplicative group). Let $\Zz_p$ denote the ring of $p$-adic integers. By \cite[Theorem 20.2.13(i)]{Hazew}, we have the following ring isomorphism:
$$\mathrm{End}_{\mathrm{formal}}(\gm,\gm)\cong \Zz_p.$$
We obtain an interesting positive characteristic version (see Example 4.15(3) and Theorem 4.16 in \cite{K8}) of raising to powers Ax-Schanuel (see Theorem \ref{risingax}). For $x\in 1+XC\llbracket X\rrbracket$ and $\alpha\in \Zz_p$ ($\mathrm{char}(C)=p>0$), we represent $\alpha$ as $\sum_{i=0}^{\infty}\alpha_ip^i$ for some $\alpha_i\in \{0,1,\ldots,p-1\}$ and we have
$$x^{\alpha}:=\lim_n \prod_{i=0}^nx^{\alpha_ip^i}.$$
\begin{theorem}\label{risingaxp}
Suppose that $\alpha\in \Zz_p$ and $[\Qq(\alpha):\Qq]>n$, Let $x_1,\dots ,x_n\in 1+XC\llbracket X\rrbracket$ be multiplicatively independent. Then
$$\trd_{C(X)}\left(x_1,\dots ,x_n,x_1^{\alpha},\dots ,x_n^{\alpha}\right)\geqslant n.$$
\end{theorem}
There is a general set-up including the additive and multiplicative cases, which we describe below following \cite{K8}.
Let us fix a positive integer $n$ and a one-dimensional algebraic group $H$ over $C$.
We introduce the following notation from \cite{K8}.
\begin{itemize}
\item Let $\mathbf{R}:=\mathrm{End}_{\mathrm{algebraic}}(H)$ and $\mathbf{S}:=\mathrm{End}_{\mathrm{formal}}(H)$.

\item We restrict our attention to algebraic groups $H$ such that $\mathbf{S}$ is a commutative domain. We regard $\mathbf{R}$ as a subring of $\mathbf{S}$.

\item Let $\mathbf{K}$ denote the field of fractions of $\mathbf{R}$ and $\mathbf{L}$ be the field of fractions of $\mathbf{S}$. We regard $\mathbf{K}$ as a subfield of $\mathbf{L}$.
\end{itemize}
\begin{example}\label{ax2ex}
In the case of the characteristic $0$, we always have $\mathbf{S}=C$, so the commutativity assumption is satisfied and we can consider any one-dimensional algebraic group as $H$.  We give some examples below.
\begin{enumerate}
\item If $H=\ga$ and characteristic is $0$, then $\mathbf{R}=\mathbf{S}=C$.

\item If $H=\ga$ and characteristic is $p>0$, then $\mathbf{R}=C[\fr]$ and $\mathbf{S}=C\llbracket \fr\rrbracket$ (see the notation introduced in the beginning of this section). This is why we needed to take $C=\Ff_p$ to ensure that $\mathbf{S}$ is commutative.

\item If $H=\gm$ and characteristic is $0$, then $\mathbf{R}=\Zz$. In the case of characteristic $p>0$, we have $\mathbf{S}=\Zz_p$ as it was mentioned above.

\end{enumerate}
\end{example}
\noindent
Below is our transcendental statement about formal endomorphisms (see \cite[Theorem 4.16.]{K8}). We need to introduce the following notions from \cite{K8}. Let $A$ be a commutative algebraic group over the field $C$ of characteristic $p>0$.
\begin{itemize}
  \item A formal map into $\widehat{A}$ is an \emph{$A$-limit map}, if it can be ``strongly approximated'' by a sequence of polynomial maps $(f_n)_n$ in the sense that the differences  $f_{n+1}-f_n$ are in the image of the $n$-th power of the appropriate Frobenius map. For example, any formal endomorphism of $\widehat{\ga}$ is a $\ga$-limit map (approximated by additive polynomials) and any formal endomorphism of $\widehat{\gm}$ is a $\gm$-limit map (approximated by multiplicative polynomials appearing in the description of $x^{\alpha}$ before the statement of Theorem \ref{risingaxp}).

\item We fix a complete $C$-algebra $\mathcal{R}$ with the residue field $C$ such that $\mathcal{R}$ is linearly disjoint from $C^{\alg}$ over $C$ and in the case of characteristic $p$ such that $L^{p^{\infty}}=C$, where $L$ is the fraction field of $\mathcal{R}$ (e.g. $\mathcal{R}$ may be the power series algebra).

  \item For $x\in A(\mathcal{R})$, we call $x$ \emph{subgroup independent}, if for any proper algebraic subgroup $A_0<A$ defined over $C$, we have $x\notin A_0(\mathcal{R})$.

\item The formal locus of $x\in A(\mathcal{R})$ over $C$ is defined as the formal subscheme of $\widehat{A}$ corresponding to the image of the map $\widehat{\mathcal{O}}_{A,0}\to \mathcal{R}$.

\item The number $\mathrm{andeg}(x)$ denotes the dimension of the formal locus of $x$ over $C$.

\end{itemize}
\begin{theorem}\label{endoax}
Take $\gamma\in \mathbf{S}$ such that $[\mathbf{K}[\gamma]:\mathbf{K}]>n$ and $\gamma:\widehat{H}\to \widehat{H}$ is an $H$-limit map. Let $\mathcal{E}:\widehat{A}\to \widehat{A}$ be the $n$-th cartesian power of $\gamma$, where $A=H^n$. Then for any subgroup independent $x\in A(\mathcal{R})_*$ we have
$$\trd_C(x,\mathcal{E}_K(x))\geqslant n+\mathrm{andeg}_C(x).$$
\end{theorem}
We showed in \cite{K8} that an unproved version of Theorem \ref{short} without the ``goodness'' assumptions implies the following conjecture. This conjecture is important for the following reasons.
\begin{itemize}
  \item If the field $C$ has characteristic $0$, then this conjecture is a theorem of Ax (\cite[Theorem 1F]{ax72}).

  \item Ax showed in \cite[Section 3]{ax72} that, in the case of characteristic $0$ (Ax did not consider the positive characteristic case), \cite[Theorem 1F]{ax72} implies  the Ax-Schanuel statements regarding the differential equation of the ``appropriate'' formal/analytic homomorphisms between algebraic groups (Ax focused on the exponential maps on semi-abelian varieties). The corresponding implication holds in the positive characteristic case as well.
\end{itemize}
\smallskip

{\bf Main Conjecture (arbitrary characteristic)}
\\
\textit{Let $G$ be an algebraic group over a field $C$ of arbitrary characteristic, $\widehat{G}$ the formalization of $G$ at the
origin and $\mathcal{A}$ a formal subgroup of $\widehat{G}$. Let $\mathcal{K}$ be a formal subscheme of $\mathcal{A}$ and let $V$ be
the Zariski closure of $\mathcal{K}$ in $G$. Then there is a formal subgroup $\mathcal{B}$ of $\widehat{G}$ which contains $\mathcal{A}$ and $\widehat{V}$ such that
$$\dim(\mathcal{B}) \leqslant  \dim(V ) + \dim(\mathcal{A}) - \dim(\mathcal{K}).$$
}

We formulate below a specific statement which would follow from the Main Conjecture above.

\smallskip

{\bf Specific Conjecture (arbitrary characteristic)}
\\
\textit{Suppose that $\mathrm{char}(C)=p>0$ and let
$$\gamma:\widehat{\gm}\to \widehat{E}$$
be a formal isomorphism where $E$ is an ordinary elliptic curve. Let $x_1,\dots ,x_n\in 1+XC\llbracket X\rrbracket$ be multiplicatively independent. Then
$$\trd_{C(X)}\left(x_1,\dots ,x_n,\gamma(x_1),\dots ,\gamma(x_n)\right)\geqslant n.$$
}

This case seems to be related to the ``interesting research paths (1)--(4)'' from Section \ref{22}. More precisely, the formal map appearing in the Specific Conjecture looks ``closest'' to the exponential map from the original Ax's Theorem, which was used by Zilber to show Weak CIT (see Section \ref{22}).
\\
\\
We finish this section with a brief discussion of the case of the Drinfeld modules. Drinfeld introduced in \cite{dr1} elliptic modules, which are now called \emph{Drinfeld modules}. Drinfeld modules can be understood as certain homomorphisms between $\Ff_q[X]$ and $K[\fr]$, where $q$ is a power of $p$ and $K=\Ff_q((\theta))$ is the non-Archimedean field of Laurent series over $\Ff_q$. An additive power series over $K$ is associated to each Drinfeld module and this series is entire on $K$. A number of transcendence results for such additive power series was obtained, see e.g. \cite{yu}. To the best of my knowledge, such results never include a version of the full Ax-Schanuel statement. For a survey of this theory, we refer the reader to \cite{Brow}.  Before the invention of Drinfeld's modules, a special case of such a series was introduced by Carlitz, which is called now the \emph{Carlitz exponential} and has the following form:
$$\exp_C=X+\sum_{i=1}^{\infty}\frac{X^{p^i}}{(\theta^{p^i}-\theta)(\theta^{p^i}-\theta^p)\ldots(\theta^{p^i}-\theta^{p^{i-1}})}.$$
 There are several Schanuel type results for the Carlitz exponential (see \cite{denis}) and a Carlitz exponential version of the (still open) conjecture on algebraic independence of logarithms of algebraic numbers was proved in \cite[1.2.6]{papa}. The power series we consider do not fit to the Drinfeld's module framework, since we consider power series with constant coefficients, that is, there is no transcendental element $\theta$ in the coefficients of our additive power series.

\section{Recent ideas and speculations}\label{secspecul}
In this section, we describe some recent early stage developments concerning Ax-Schanuel type problems. One of them regards combining the results from \cite{cfn3} with Ax's theorem on the dimension of intersections. The other one is about differential versions of Ax-Schanuel in positive characteristic.

\subsection{Towards a general statement of Ax-Schanuel}\label{specul1}
Ax-Schanuel statements for analytic/formal homomorphisms in the case of characteristic $0$ have one ``umbrella statement'' from which they all follow, which is Ax's theorem on the dimension of intersections from Section \ref{seccharp}. No such an ``umbrella statement'' was known for  Ax-Schanuel statements for the maps like the $j$-invariant map, until the recent preprint \cite{cfn3}, where a general form of an Ax-Schanuel type result is given (\cite[Theorem A]{cfn3}). In this statement, the algebraic group $G$ is again back in the picture (e.g. $G=\mathrm{PGL}_2(\Cc)$ in the case of the $j$-invariant map), but the statement is quite technical and it is phrased in terms of leaves of flat connections on $G$-principal bundles, where such a leaf plays a role of the analytic subgroup $\mathcal{A}$ from Ax's theorem on the dimension of intersections from Section \ref{seccharp}.

\smallskip

{\bf Connection version of Ax-Schanuel (\cite[Theorem A]{cfn3})}
\\
\textit{Let $\nabla$ be a $G$-principal flat connection on the algebraic bundle $P\to Y$ such that:
\begin{itemize}
  \item the algebraic group $G$ is sparse;
  \item the Galois group of $\nabla$ coincides with $G$.
\end{itemize}
Let $V$ be an algebraic subvariety of $P$ and $\mathcal{L}$ be a horizontal leaf of $\nabla$. If
$$\dim V < \dim(V \cap \mathcal{L})+\dim G$$
then the projection of $V\cap \mathcal{L}$ in $Y$ is contained in a $\nabla$-special subvariety of $Y$.}

\smallskip

Sparsity of the algebraic group $G$ above means that there are no proper Zariski dense complex analytic subgroups of $G$. The notion of a ``$\nabla$-special'' is more technical, it is phrased in terms of the Galois group of a connection (see \cite[Definition 2.4]{cfn3}).

Unlike in the case of \cite[Theorem 1]{ax72}, no analytic subgroup appears in \cite[Theorem A]{cfn3}, so this theorem does not generalize \cite[Theorem 1]{ax72}. We propose such a generalization which encompasses both the Connection version of Ax-Schanuel and \cite[Theorem 1]{ax72}. It will appear in \cite{goko}.

\smallskip

{\bf Connection and subgroup Ax-Schanuel}
\\
\textit{Let $\nabla$ be a $G$-principal flat connection on the algebraic bundle $P\to Y$ such that the Galois group of $\nabla$ coincides with $G$ and
\begin{itemize}
  \item $V$ be an algebraic subvariety of $P$,
  \item   $\mathcal{A}$ be an analytic subgroup of $G$,
  \item $\mathcal{L}$ be a horizontal leaf of $\nabla$.
\end{itemize}
Suppose that $\mathcal{V}$ is an analytic submanifold of $\mathcal{A}$ which is Zariski dense in $V$. If
$$\dim V < \dim(V \cap \mathcal{L})+\dim G$$
then there is an analytic subgroup $\mathcal{H}$ of $G$ such that
$$\dim \mathcal{H} < \dim(V)-\dim(\mathcal{V})$$
and $V\subseteq \mathcal{A}\mathcal{H}$.}
\smallskip

The results mentioned above concern the case of characteristic zero. In the ``Main Conjecture'' from Section \ref{seccharp}, the notion of ``analytic'' is replaced with the notion of ``formal'' (see Remark \ref{lookpos}) which makes sense in the case of arbitrary characteristic. The connection version of Ax-Schanuel (\cite[Theorem A]{cfn3}) mentioned above has not been considered in the positive characteristic case before, since it requires an appropriate version of the notion of a connection in positive characteristic. This is work in progress (\cite{goko}).

\subsection{Hasse-Schmidt differential Ax-Schanuel}
Positive characteristic versions of the differential Ax's theorem have not been studied yet. It is clear that we can not consider the usual derivations anymore, since the constants of differential fields of positive characteristic contain the image of the Frobenius map, hence there is no room for any transcendence. It looks natural in this case to replace the derivations with iterative Hasse-Schmidt derivations and the field of constants with the field of absolute constants, we give the necessary definitions below.
\begin{itemize}
  \item A sequence $\partial = (\partial_{n}:R\to R)_{n\in \Nn}$ of additive maps on a ring is
called an \emph{HS-derivation} if $\partial_{0}$ is the identity map, and for all $n\in \Nn$ and $x,y\in R$, we have:
$$\partial_{n}(xy) = \sum_{i+j=n}\partial_{i}(x)\partial_{j}(y).$$

  \item An HS-derivation $\partial$ is called \emph{iterative} if for all $i,j\in \Nn$ we have
$$\partial_{i}\circ\partial_{j} = {i+j \choose i}\partial_{i+j}.$$

  \item If $(K,\partial)$ is a field with a Hasse-Schmidt derivation, then its field of absolute constants is the following intersection:
  $$\bigcap_{i=1}^{\infty}\ker\left(\partial_i\right).$$
\end{itemize}
The passages between the differential Ax-Schanuel and the power series Ax-Schanuel (described in Remark \ref{passege}) work only one way for positive characteristic case, since the power series ring has a natural iterative Hasse-Schmidt derivation on it. However, it is not clear how to proceed in the opposite way, so Hasse-Schmidt differential Ax-Schanuel type results need to be proved separately. This is work in progress (\cite{goko}).

We state below two such results which will appear in \cite{goko} to give a flavour of these kind of Ax-Schanuel conditions. Assume that $(K,\partial)$ is a field of characteristic $p>0$ with a Hasse-Schmidt derivation and $C$ is a field contained in the field of absolute constants of $(K,\partial)$.
\\
\\
\textbf{Additive Hasse-Schmidt differential Ax-Schanuel}
\\
\\
\textit{Let
$$F=\sum_{m=0}^{\infty} c_mX^{p^m}\in \Ff_p\llbracket \fr\rrbracket$$
and assume that the algebraic degree of $F$ over $\Ff_p(\fr)$ is greater than $n$. Take $x_1,\ldots,x_n,y_1,\ldots,y_n\in K$ such that $x_1,\ldots,x_n$ are linearly independent over $\Ff_p[\fr]$ and for all $i\in \{1,\ldots,n\}$:
$$D_1\left(y_i-c_0x_i\right)=0,$$
$$D_p\left(y_i-c_0x_i-c_1x_i^p\right)=0,$$
$$\ldots$$
$$D_{p^m}\left(y_i-c_0x_i-c_1x_i^p-\ldots -c_mx^{p^m}\right)=0,$$
$$\ldots$$
Then we have:
$$\td_{\Ff_p}(x_1,y_1\ldots,x_n,y_n)\geqslant n+1.$$}
\mbox{}
\\
\textbf{Multiplicative Hasse-Schmidt differential Ax-Schanuel}
\\
\\
\textit{Let
$$\gamma=\sum c_ip^i\in \Zz_p$$
and assume that the algebraic degree of $\gamma$ over $\Qq$ is greater than $n$.
Take $x_1,\ldots,x_n,y_1,\ldots,y_n\in K$ such that $x_1,\ldots,x_n$ are multiplicatively independent and for all $i\in \{1,\ldots,n\}$:
$$D_1\left(y_ix_i^{-c_0}\right)=0,$$
$$D_p\left(y_ix_i^{-c_0-c_1p}\right)=0,$$
$$\ldots$$
$$D_{p^m}\left(y_ix_i^{-c_0-c_1p-\ldots -c_m{p^m}}\right)=0,$$
$$\ldots$$
Then we have:
$$\td_{C}(x_1,y_1\ldots,x_n,y_n)\geqslant n+1.$$}

\bibliographystyle{plain}
\bibliography{harvard}

\end{document}